\documentclass[twoside]{article}
\input{spsp10.sty}
\begin{document}
\thispagestyle{empty}
\fancyhead[LE]{Goodarzi, F. and Zamini, R.}
\fancyhead[RO]{Estimation of varextropy}
\begin{center}
\vspace*{4cm}
{\Large \bf  On the estimation of varextropy under complete data }\\
{\bf
F. Goodarzi$^{a}$,
R. Zamini$^{b}$\footnote{{\bf Raheleh Zamini:}{\tt  zamini@khu.ac.ir}}\\
$^{a}$Faculty of Mathematical Sciences, Department of Statistics, University of Kashan, Kashan, Iran.\\
$^{b}$Department of Mathematics, Faculty of Mathematical Sciences and Computer, Kharazmi University, Tehran, Iran.
}
\end{center}
\noindent{\bf Abstract:} In this paper, we propose nonparametric estimators for varextropy function of an absolutely continuous random variable. Consistency  
of the estimators is established under suitable regularity conditions. Moreover, a  simulation study is performed to compare the performance of the proposed estimators based on mean squared error (MSE) and bias. Furthermore, by using the proposed estimators some tests are constructed for uniformity. 
It is shown that the varextropy-based test proposed here performs well compared to the power of the 
other uniformity hypothesis tests.


\section{Introduction}
Suppose that $X$ be a non-negative and an absolutely continuous random variable with the probability density function (pdf) $f$, distribution function $F$ and survival function $\overline F$. A prevalent uncertainty measure is the expectation of $-\log f(X)$, as defined by \cite{Shannon}, and  given by
\begin{align}
H(X)=-\int_{0}^{+\infty}f(x)\log f(x)dx,
\end{align}
Many researchers have estimated the
entropy of continuous random variables in their studies. 
Some of these researchers are \cite{Vasicek}, \cite{Van}, \cite{Ebrahimi}, \cite{Correa} and \cite{Zamanzade}.
 The varentropy of a random variable $X$  is defined as 
\begin{align}\label{VEn}
VE(X)&=Var[-\log f(X)]
\nonumber\\&=
\int_{0}^{+\infty}f(x)[\log f(x)]^2dx-\left[\int_{0}^{+\infty}f(x)\log f(x)dx\right]^2.
\end{align} 
Hence the variability in the information content of $X$ is  measured by the varentropy. The relevance
of this measure to physics, computer and mathematics sciences has been pointed out in various studies, including \cite{Madiman}, \cite{Jiang}, \cite{Fradelizi} and \cite{Li}. 
Recently  \cite{Alizadeh} proposed some estimators for varentropy of a continuous random variable. 

The notion of entropy is recently intertwined with a complementary dual measure, designated as extropy, by \cite{Lad}. 
The extropy of the random variable $X$ is defined as:
\begin{align}
J(X)
=-\frac{1}{2}\int_{0}^{+\infty}f^2(x)dx.
\end{align}
It is clear that $J(X)\leq 0$.  Scoring the forecasting distributions using the total log-scoring rule is one of the statistical applications of extropy (see \cite{Gne}).
Furthermore,  the extropy has been universally investigated in commercial or scientific areas such as astronomical measurements of heat distributions in galaxies,  (see \cite{Furu} and \cite{Von}).
The concept of extropy is useful in automatic speech recognition (Refer to \cite{Bec}). 
Moreover, extropy is a measure better than entropy in some scenarios in statistical mechanics and thermodynamics.
For more studies on extropy, see \cite{Qiu2018b}, \cite{Yang}.
Recently, the problem of estimating $J(X)$ has been noted by \cite{Qiu2018a}, \cite{AliJar} and \cite{al}. 

\cite{Vasel} defined the varextropy for the absolutely continuous random variable $X$ as: 
\begin{align}\label{vex}
VJ(X)=Var\left[-\frac{1}{2}f(X)\right]=\frac{1}{4}E[f^2(X)]-J^2(X).
\end{align}
Moreover, \cite{Good} obtained lower bounds for varextropy. The varextropy measure can be used as an alternative 
measure to Shannon entropy and a measure of the amount of uncertainty.

The rest of the work is structured as follows. In Section 2, we propose nonparametric estimators for varextropy function and investigate consistency of the estimators under suitable conditions.
We present the results of the Monte Carlo studies on biases and mean squared errors (MSEs) of the proposed estimators in Section 3.
In Section 4, we introduce some goodness-of–fit tests of uniformity by using the our proposed estimators and
 compare the powers of our proposed test statistics with some known test statistics.
We apply the newly introduced test procedures  to real data sets for 
illustration in Section 5.

\section{Varextropy estimation for complete data}
In this section, we propose five estimators for varextropy function under complete data and investigate some of their asymptotic properties.
\subsection{The first estimator}
Before constructing the  first estimator for $VJ(X)$, we first review some required concepts and notations. Let $X$ be a continuous random variable 
with density function $f$ and distribution function $F$. Consider the class of functions based on $f(x)$ as 
\begin{align}
T(f)=\int_{-\infty}^{+\infty}f(x)\Phi(f(x))w(x)dx,
\end{align}
where $\Phi$ is a real valued differentiable function on $[0, \infty)$ and $w$ is a real valued function on $[0, \infty$).

Assuming that $f$ is bounded, \cite{Van} proposed the following estimator for $T(f)$:
\begin{align}\label{eq2}
T_{m, n}=\frac{1}{2(n-m)}\sum_{j=1}^{n-m}\Phi\left(\frac{m}{n+1}(X_{(j+m)}-X_{(j)})^{-1}\right)\left(w(X_{(j)})+w(X_{(j+m)})\right),
\end{align}
where $X_{(1)}\leq \cdots\leq X_{(n)}$ are the order statistics corresponding to the sample $X_1, \ldots, X_n$ and $m$ is an integer such that $1\leq m\leq n$. The estimator given in \eqref{eq2} is constructed according to the fact that $f(X_{(j)})$, the density at the point $X_{(j)}$, is replaced by $m\left((n+1)(X_{(j+m)}-X_{(j)})\right)^{-1}$, where is the local histogram estimate of $f(X_{(j)})$. 

Under some conditions, \cite{Van} showed that $T_{m, n}$ is a.s. consistent, i.e., 
\begin{align}
T_{m, n}\to T(f),\hspace{0.5cm}a.s.\hspace{0.5cm}as\hspace{0.5cm}m, n\to\infty.
\end{align}\label{eq3}
(See Theorem 1 of \cite{Van} for some details.) Two special cases of $T(f)$ are 
\begin{align}\label{eq4}
T^{\prime}(f)=\int_{0}^{+\infty}f^3(x)dx, 
\end{align}
and 
\begin{align}\label{eq5}
T^{\prime\prime}(f)=\int_{0}^{+\infty}f^2(x)dx, 
\end{align}
obtained by taking $\Phi(x)$ equal to $x^2$ and $x$, respectively, and $w$  is identically equal to $I_{[0, \infty)}(x)$. 
According to $T_{m, n}$ defined in \eqref{eq2}, estimators for $T^{\prime}(f)$ and $T^{\prime\prime}(f)$ 
can be proposed as 
\begin{align}\label{eq6}
T^{\prime}_{m, n}=\frac{1}{n-m}\sum_{j=1}^{n-m}{\Big(\frac{m}{n+1}{(X_{(j+m)}-X_{(j)})}^{-1}\Big)}^2,
\end{align}
and 
\begin{align}\label{eq7}
T^{\prime\prime}_{m, n}=\frac{1}{n-m}\sum_{j=1}^{n-m}{\Big(\frac{m}{n+1}{(X_{(j+m)}-X_{(j)})}^{-1}\Big)},
\end{align}
respectively.

Now, the expression given in  \eqref{vex}   for $VJ(X)$, motivates us to construct the first estimator of varextropy based on $T^{\prime}_{m, n}$ and  $T^{\prime\prime}_{m, n}$ as follows:
\begin{align}\label{eq8}
VJV_{m, n}=\frac{1}{4}\{T^{\prime}_{m, n}-{(T^{\prime\prime}_{m, n})}^2\}.
\end{align}
In the following theorem, we establish the a.s. consistency of $VJV_{m, n}$.
\begin{thm}
Let the density function $f$ be bounded, then $VJV_{m, n}\to VJ(X) \hspace{0.2cm}a.s.,$ provided $m, n\to\infty$, $\frac{m}{\log n}\to\infty$ and $\frac{m}{n}\to 0$.
\end{thm}
\begin{proof}
\begin{align}\label{eq9}
VJV_{m,n}-VJ(X)&=\frac{1}{4}\{T^{\prime}_{m, n}-{(T^{\prime\prime}_{m, n})}^2\}-
\frac{1}{4}\{T^{\prime}(f)-{(T^{\prime\prime}(f))}^2\}\nonumber
\\&=\frac{1}{4}\{T^{\prime}_{m, n}-T^{\prime}(f)\}-
\frac{1}{4}\{{(T^{\prime\prime}_{m, n})}^2-{(T^{\prime\prime}(f))}^2\}.
\end{align}
Theorem 1 of \cite{Van} ensures that 
\begin{align}\label{eq10}
T^{\prime}_{m, n}\to T^{\prime}(f)\hspace{0.2cm}a.s.\hspace{0.2cm}as\hspace{0.2cm}m, n\to\infty,
\end{align}
and
\begin{align}\label{eq11}
(T^{\prime\prime}_{m, n})^2\to (T^{\prime\prime}(f))^2\hspace{0.2cm}a.s.\hspace{0.2cm}as\hspace{0.2cm} m, n\to\infty.
\end{align}
\eqref{eq9}, \eqref{eq10} and \eqref{eq11} complete the proof.
\end{proof}
\subsection{The second estimator}
Let $\hat f(x)$ be the kernel density function of $f$ defined by
\begin{align}\label{eq12}
\hat f(x)=\frac{1}{nh_n}\sum_{i=1}^nK\Big(\frac{x-X_i}{h_n}\Big),
\end{align}
where $K$ is the kernel density function and $h_n$ is the bandwidth parameter. The second estimator for $VJ(X)$ can be proposed based on $\hat f(x)$ as:
\begin{align}\label{eq13}
VJD=\frac{1}{4}\left[\int{\hat f}^3 (x)dx-\Big(\int{\hat f}^2(x)dx\Big)^2\right].
\end{align}
The following assumptions on $K$ are used in the Theorem \ref{Th2} of this section.

$(A_1):$ $K$ is uniformly continuous and of bounded variation $V(K)$.

$(A_2): K(x)\to 0\hspace{0.2cm}as\hspace{0.2cm}|x|\to\infty$. 

$(A_3): \int |x \log |x||^{\frac{1}{2}}|dK(x)|<\infty$.
\begin{thm}\label{Th2}
Suppose $K$ satisfies Assumptions $A_1-A_3$ and $f$ is uniformly continuous. If $h_n\to 0$ and ${(nh_n)}^{-1}\log n\to0$ as $n\to\infty$, then 
\begin{align*}
VJD\to VJ(X)\hspace{0.2cm}a.s.\hspace{0.2cm}as \hspace{0.2cm}n\to\infty.
\end{align*}
\end{thm}
\begin{proof}
From \eqref{vex} and \eqref{eq13} we have
\begin{align}\label{eq14}
VJD-VJ(X)=I_n+II_n, 
\end{align}
where 
\begin{align}\label{eq15}
I_n=\frac{1}{4}\left(\int_{0}^{+\infty}{\hat f}^3(x)dx-\int_{0}^{+\infty}{f^3(x)}dx\right),
\end{align}
and 
\begin{align}\label{eq16}
II_n=\left(-\frac{1}{2}\int_{0}^{+\infty}f^2(x)dx\right)^2-\left(-\frac{1}{2}\int_{0}^{+\infty}{\hat f}^2(x)dx\right)^2.
\end{align}
Observe that
\begin{align}\label{eq17}
|I_n|&=\frac{1}{4}\left|\int_{0}^{+\infty}{\hat f}^3(x)dx-\int_{0}^{+\infty}f^3(x)dx\right|\nonumber\\
&\leq \frac{1}{4}\left|\int_{0}^{+\infty}{\hat f}^2(x)\Big(\hat f(x)-f(x)\Big)dx\right|\nonumber
\\&+\frac{1}{4}
\left|\int_{0}^{+\infty}\Big({\hat f}^2(x)-f^2(x)\Big)f(x)dx\right|.
\end{align}
Under Assumptions  $A_1-A_{3}$, from Theorem A of \cite{Silv}, one can see that 
\begin{align}\label{eq18}
\lim_{n\to\infty}\sup_{x}|\hat f(x)-f(x)|=0, \ \ a.s.
\end{align}
On the other hand, a simple algebra calculation shows that 
\begin{align}\label{eq19}
{\hat f}^2(x)-f^2(x)=\sum_{i=1}^{2}{(f(x))}^{i-1}{(\hat f(x))}^{2-i}(\hat f(x)-f(x)).
\end{align}
\eqref{eq18} and \eqref{eq19} ensure that 
\begin{align}\label{eq20}
\lim_{n\to\infty}\sup_{x}|{\hat f}^2(x)-f^2(x)|=0, \ \ a.s.
\end{align}
From \eqref{eq17}, \eqref{eq18} and \eqref{eq20}, one can observe that
\begin{align}\label{eq21}
\lim_{n\to \infty}I_{n}=0,\hspace{0.5cm}a.s.
\end{align}
Next, to deal with $II_n$, observe that 
\begin{align}\label{eq22}
\left|\int_{0}^{+\infty}\hat f^2(x)dx-\int_{0}^{+\infty} f^2(x)dx\right|&\leq \left|\int_{0}^{+\infty}\hat f^2(x)dx-\int_{0}^{+\infty}\hat f(x)dF(x)\right|\nonumber\\&+\left|\int_{0}^{+\infty}\hat f(x)dF(x)-\int_{0}^{+\infty}f^2(x)dx\right|\nonumber
\\&=\left|\int_{0}^{+\infty}\hat f(x)(\hat f(x)-f(x))dx\right|\nonumber\\&+\left|\int_{0}^{+\infty}(\hat f(x)-f(x))dF(x)\right|.
\end{align}
Since from \eqref{eq18}, $\hat f(x)$ is a.s. consistent, we can write
\begin{align}\label{eq23}
\lim_{n\to \infty} \int_{0}^{+\infty}\hat f^2(x)dx=\int_{0}^{+\infty}f^2(x)dx,\hspace{0.2cm}a.s.
\end{align}
As a result 
\begin{align}\label{eq24}
\lim_{n\to \infty}II_n=0, \hspace{0.2cm}a.s.
\end{align}
\eqref{eq14}, \eqref{eq21} and \eqref{eq24} complete the proof.
\end{proof}
\subsection{The third estimator}
The third estimator is proposed by:
\begin{align}\label{eq26}
\frac{1}{4}\left(\frac{1}{n}\sum_{i=1}^nf^2(X_i)-\Big(\frac{1}{n}\sum_{i=1}^nf(X_i)\Big)^2\right),
\end{align}
using the empirical distribution function $F_n(x)=\frac{1}{n}\sum_{i=1}^nI(X_i\leq x)$.

Since $f$ is unknown, a natural estimator for $VJ(X)$ can be obtained by substituting $f$ by $\hat f$ defined in \eqref{eq12}, and is given by 
\begin{align}\label{eq27}
VJB_n&=\frac{1}{4}\left(\int_{0}^{+\infty}\hat f^2(x)dF_n(x)-\Big(\int_{0}^{+\infty}\hat f(x)dF_n(x)\Big)^2\right)
\nonumber\\&=\frac{1}{4}\left(\frac{1}{n}\sum_{i=1}^{n}(\tilde f_n(X_i))^2-\Big(\frac{1}{n}\sum_{i=1}^{n}\tilde f_n(X_i)\Big)^2\right),
\end{align}
where 
\begin{align}\label{eq28}
\tilde f_n(X_i)=\frac{1}{(n-1)h_{n-1}}\sum_{j=1, j\ne i}^nK\left(\frac{X_i-X_j}{h_{n-1}}\right).
\end{align} 
In the following theorem, we investigate a.s. consistency of $VJ_n(X)$. 
\begin{thm}
Suppose that $\lim_{x\to\infty}f(x)$ exists and is finite. If $f$ is bounded, then under the conditions of Theorem \ref{Th2}, we can write:
\begin{align*}
\lim_{n\to \infty}|VJB_n-VJ(X)|=0,\hspace{0.2cm}a.s.
\end{align*} 
\end{thm}
\begin{proof}
By usng \eqref{vex} and \eqref{eq27} we have
\begin{align}\label{eq29}
|VJB_n-VJ(X)|\leq |I^{\prime}_n|+|I^{\prime\prime}_n|,
\end{align}
where 
\begin{align}\label{eq30}
I^{\prime}_n=\frac{1}{4n}\sum_{i=1}^{n}(\tilde f_n(X_i))^2-\frac{1}{4}\int_{0}^{+\infty}f^2(x)dF(x),
\end{align}
and 
\begin{align}\label{eq31}
I^{\prime\prime}_n=\left(-\frac{1}{2}\int_{0}^{+\infty}f^2(x)dx\right)^2-\left(-\frac{1}{2n}\sum_{i=1}^{n}(\tilde f_n(X_i))\right)^2.
\end{align}
Using the definition of step function $F_n(x)$ and adding and subtracting $\frac{1}{4}\int_{0}^{+\infty}f^2(x)dF_n(x)$, $I_n^{\prime}$ can be rewritten as 
\begin{align}\label{eq321}
I^{\prime}_n&=\frac{1}{4}\int_{0}^{+\infty}\hat f^2(x)dF_n(x)-\frac{1}{4}\int_{0}^{+\infty}f^2(x)dF(x)\nonumber
\\&=\frac{1}{4}\int_{0}^{+\infty}(\hat f^2(x)-f^2(x))dF_n(x)+\frac{1}{4}\int_{0}^{+\infty}f^2(x)d(F_n(x)-F(x)).
\end{align}
On the other hand, since $f$ is bounded and $\lim_{x\to\infty}f(x)$ exists and is finite, by using the integration by parts we have:
\begin{align}\label{eq322}
\frac{1}{4}\int_{0}^{+\infty}f^2(x)d(F_n(x)-F(x))=-\frac{1}{2}\int_{0}^{+\infty}(F_n(x)-F(x))f(x)f^{\prime}(x)dx.
\end{align}
From \eqref{eq321} and \eqref{eq322}, we can write:
\begin{align}\label{eq323}
I^{\prime}_n=\frac{1}{4}\int_{0}^{+\infty}(\hat f^2(x)-f^2(x))dF_n(x)-\frac{1}{2}\int_{0}^{+\infty}(F_n(x)-F(x))f(x)f^{\prime}(x)dx. 
\end{align}
Hence
\begin{align}\label{eq324}
|I_n^{\prime}|\leq \frac{1}{4}\sup_{x}|\hat f^2(x)-f^2(x)|+\frac{1}{4}\sup_{x}|F_n(x)-F(x)|(\lim_{x\to\infty}f^2(x)+(\sup_x f(x))^2).
\end{align}
 
The a.s. convergence of $\hat f^2(x)$ given in \eqref{eq20} and a.s. convergence of $F_n(x)$ given in the Glivenko-Cantelli theorem ensure that 
\begin{align}\label{eq33}
\lim_{n\to \infty}|I^{\prime}_n|=0.\hspace{0.2cm}
\end{align}
In a similar way, using the triangle inequality and adding and subtracting $\int_0^{+\infty} f(x)dF_n(x)$, we can see that 
\begin{align}\label{eq34}
\left|\int_{0}^{+\infty}\hat f(x)dF_n(x)-\int_{0}^{+\infty} f(x)dF(x)\right|&\leq 
\left|\int_{0}^{+\infty}(\hat f(x)-f(x))dF_n(x)\right|\nonumber\\&+ \left|\int_{0}^{+\infty} f(x)d(F_n(x)-F(x))
\right|.
\end{align}
Since $f$ is bounded and $\lim_{x\to\infty}f(x)$ exists and is finite, using the integration by parts we can write:
\begin{align}\label{eq341}
\int_{0}^{+\infty} f(x)d(F_n(x)-F(x))=-\int_{0}^{+\infty}(F_n(x)-F(x))f^{\prime}(x)dx.
\end{align}
 \eqref{eq34} and \eqref{eq341} imply that 
\begin{align}\label{eq342} 
\left|\int_{0}^{+\infty}\hat f(x)dF_n(x)-\int_{0}^{+\infty} f(x)dF(x)\right|&\leq  
\sup_x|\hat f(x)-f(x)|\nonumber\\&+\sup_x|F_n(x)-F(x)|(\lim_{x\to\infty}f(x)+\sup_x f(x)).
\end{align}
It can be concluded from the a.s. convergence of $\hat f$ given in \eqref{eq18}, the Glivenko-Cantelli theorem and the finitness of $\lim_{x\to\infty}f(x)$ and $\sup_xf(x)$ that 
\begin{eqnarray}\label{eq35}
\lim_{n\to\infty}|I^{\prime\prime}_n|=0.
\end{eqnarray}
\eqref{eq29},  \eqref{eq33} and \eqref{eq35} complete the proof.
\end{proof}
\subsection{The fourth estimator}
Let $Q(u)=\inf\{x: F(x)\geq u\}, 0\leq u\leq 1,$ be the quantile function corresponding to the distribution function $F(x)$ and 
$q(u)=Q^{\prime}(u)=\frac{1}{f(Q(u))}$ be its quantile density function.  \cite{Soni} proposed a smooth estimator for $q(u)$ given by 
\begin{align}
\tilde q_n(u)&=\frac{1}{h_n}\int_{0}^{1}\frac{K(\frac{t-u}{h_n})}{f_n(Q_n(t))}dt\nonumber\\&=\frac{1}{nh_n}\sum_{i=1}^n\frac{K(\frac{S_i-u}{h_n})}{f_n(X_{(i)})},
\end{align} 
where $S_i$ is the proportion of observations less than or equal to $X_{(i)}$, the $ith$  order statistic, $f_n(u)$ is the same $\hat f(u)$ given in \eqref{eq12}, and $Q_n(u)=\inf\{x: F_n(x)\geq u\}, 0\leq u\leq 1,$ is the empirical estimator of $Q(u)$.

The fourth estimator for varextropy function $VJ(X)$ is obtained by considering the fact that $VJ(X)$ can be rewritten  in terms of the quantile function as follows: 
\begin{align}\label{eq37}
VJ(X)=\frac{1}{4}\left[\int_{0}^{1}f^2(Q(u))du-\left(\int_{0}^{1}f(Q(u))du\right)^2\right]\nonumber\\=
\frac{1}{4}\left[\int_{0}^{1}\frac{du}{(q(u))^2}-\left(\int_{0}^{1}\frac{du}{q(u)}\right)^2\right].
\end{align}
We consider the following plug-in estimator for quantile-based varextropy, given by 
 \begin{align}\label{eq38}
VJS_n=
\frac{1}{4}\left[\int_{0}^{1}\frac{du}{(\tilde q_n(u))^2}-\left(\int_{0}^{1}\frac{du}{\tilde q_n(u)}\right)^2\right].
\end{align}
In the following theorem, we prove the a.s. consistency of $VJS_n$.
\begin{thm}
Let $f$ be bounded, then
\begin{align*}
VJS_n\to VJ(X), \ \ a.s. \ \ n\to\infty.
\end{align*}
\end{thm}
\begin{proof}
We have
\begin{align}\label{eq39}
VJS_n- VJ(X)&=\frac{1}{4}\left[\int_{0}^{1}\frac{du}{(\tilde q_n(u))^2}-\int_{0}^{1}\frac{du}{(q(u))^2}\right]
\nonumber\\&-\frac{1}{4}\left[\left(\int_{0}^{1}\frac{du}{\tilde q_n(u)}\right)^2-\left(\int_{0}^{1}\frac{du}{q(u)}\right)^2\right]\nonumber\\&:=J_1+J_2.
\end{align}
Now
\begin{align}\label{eq40}
\int_{0}^{1}\left(\frac{1}{(\tilde q_n(u))^2}-\frac{1}{(q(u))^2}\right)du=\int_{0}^{1}(q(u)-\tilde{q}_n(u)).\frac{q(u)+\tilde{q}_n(u)}{(\tilde q_n(u))^2q^2(u)}du.
\end{align}
Using proof of Theorem 3.4 of \cite{Sub}, $\sup_{t}|\tilde q_n(t)-q(t)|\to 0\ \ as \ \ n\to\infty$, hence the above right expression asymptotically reduces to
\begin{align}
\int_{0}^{1}(q(u)-\tilde q_n(u))\frac{2}{q^3(u)}du.
\end{align}  
Since $f$ is bounded, there exist $M>0$, such that $|f|\leq M$ and hence 
\begin{align}
\int_{0}^{1}(q(u)-\tilde q_n(u))\frac{2}{q^3(u)}du\leq 2\sup_u|q(u)-\tilde q_n(u)|\int_{0}^1f^3(Q(u))du\leq 2M^3\sup_u|q(u)-\tilde q_n(u)|\to 0\ \ as\ \ n\to\infty.
\end{align}     
As a result
\begin{align}\label{eq42}
J_1\to 0 \ \ as \ \ n\to\infty.
\end{align}
In a similar way, we can show that 
\begin{align}\label{eq43}
J_2\to 0 \ \ as \ \ n\to\infty.
\end{align}
\eqref{eq39}, \eqref{eq42} and \eqref{eq43}, ensure that 
\[VJS_n\to VJ(X)\ \ \ as \ \ n\to \infty.\]
\end{proof}
\subsection{The fifth estimator}
By analogy to \cite{Qiu2018a} and using the \eqref{eq37}, we propose the following estimator for varextropy function: 
\begin{align}
VJQ_{m, n}=\frac{1}{4}\left[\frac{1}{n}\sum_{i=1}^n{\left(\frac{c_i\frac{m}{n}}{X_{(i+m)}-X_{(i-m)}}\right)}^2-{\left(\frac{1}{n}\sum_{i=1}^n\frac{c_i\frac{m}{n}}{X_{(i+m)}-X_{(i-m)}}\right)}^2\right],
\end{align}
where the window size $m$ is a positive integer smaller than $n/2$, $X_{(i)}=X_{(1)}$ if $i<1$,  $X_{(i)}=X_{(n)}$ if $i>n$ and 
\[
c_i=\left\{ 
\begin{array}{ll}
1+\frac{i-1}{m},&1\leq i\leq m,\\
2,& m+1\leq i\leq n-m,\\
1+\frac{n-i}{m}, &n-m+1\leq i\leq n.
\end{array}\right.
\]
The consistency of the proposed estimator is shown in the following theorem.
\begin{thm}\label{thm2.5}
Let $X_1$, $\ldots$, $X_n$ be a random sample from distribution function $F$ with bounded probability density function $f$ and finite variance. Then we can write:
\begin{align}
VJQ_{m, n}\stackrel{\rm p}{\to}VJ(X), \hspace{0.3cm}as \ \ n\to\infty, \ \ m\to\infty, \ \ \frac{m}{n}\to 0.
\end{align}
\end{thm}
\begin{proof}
A simple algebra shows that:
\begin{align}\label{1t}
VJQ_{m, n}-VJ(X):=I_1+I_2,
\end{align}
where 
\begin{align}
I_1=\frac{1}{4}\left[\frac{1}{n}\sum_{i=1}^n{\left(\frac{c_i\frac{m}{n}}{X_{(i+m)}-X_{(i-m)}}\right)}^2-\int_{0}^{+\infty} f^3(x)dx\right],
\end{align}
and 
\begin{align}
I_2=\frac{1}{4}\left[{\left(\frac{1}{n}\sum_{i=1}^n\frac{c_i\frac{m}{n}}{X_{(i+m)}-X_{(i-m)}}\right)}^2-
\left(\int_{0}^{+\infty} f^2(x)dx\right)^2
\right].
\end{align}
Theorem 2.1 of \cite{Qiu2018a}, clearly shows that $I_2\stackrel{p}{\to}0.$ 
For $I_1$, first we write $\frac{1}{n}\sum_{i=1}^n{\Big(\frac{c_i\frac{m}{n}}{X_{(i+m)}-X_{(i-m)}}\Big)}^2$ as follows 
\begin{align}
\frac{1}{n}\sum_{i=1}^n(c_i\frac{m}{n})^2
{\Big(\frac{1}{X_{(i+m)}-X_{(i-m)}}\Big)}^2=
\frac{1}{n}\sum_{i=1}^n{\Big(\frac{\frac{2m}{n}}{X_{(i+m)}-X_{(i-m)}}\Big)}^2\nonumber\\-
\frac{2}{n}\sum_{i=1}^n\xi_{im}\frac{2m}{n}\frac{1}{X_{(i+m)}-X_{(i-m)}}(1-\frac{c_i}{2}):=J_1-J_2.
\end{align}
where $\xi_{im}$ is a random point between the points $\frac{c_i\frac{m}{n}}{X_{(i+m)}-X_{(i-m)}}$ and $\frac{\frac{2m}{n}}{X_{(i+m)}-X_{(i-m)}}$.

Next to deal with $J_1$, observe that 
\begin{align}\label{2t}
J_1&=\frac{1}{n}\sum_{i=1}^{n}{\Big(\frac{\frac{2m}{n}}{X_{(i+m)}-X_{(i-m)}}\Big)}^2=\frac{1}{n}\sum_{i=1}^{n}{\Big(\frac{F(X_{(i+m)})-F(X_{(i-m)})}
{X_{(i+m)}-X_{(i-m)}}\Big)}^2\left[{\Big(\frac{\frac{2m}{n}}{F(X_{(i+m)})-F(X_{(i-m)}}\Big)}^2-1\right]\nonumber\\&+\frac{1}{n}\sum_{i=1}^{n}{\Big(\frac{F(X_{(i+m)})-F(X_{(i-m)})}{X_{(i+m)}-X_{(i-m)})}\Big)}^2\nonumber\\&:=J_{11}+J_{12}. 
\end{align}
Now, since $f(x)$ is positive and continuous, there exists a value $X_i\in (X_{(i-m)}, X_{(i+m)})$ such that 
\begin{align}\label{3t}
f(X_i)=\frac{F(X_{(i+m)})-F(X_{(i-m)})}{X_{(i+m)}-X_{(i-m)}},
\end{align}
hence
\begin{align}\label{4t}
J_{11}=\frac{1}{n}\sum_{i=1}^{n}{(f(X_i))}^2\left[{\Big(\frac{\frac{2m}{n}}{F(X_{(i+m)})-F(X_{(i-m)})}\Big)}^2-1\right].
\end{align}
On the other hand, since $F_n(X_{(i+m)})-F_n(X_{(i-m)})=\frac{2m}{n}$, Chung's law of the iterated logarithm (LIL) implies that:
\begin{align}
\frac{F(X_{(i+m)})-F(X_{(i-m)})}{F_n(X_{(i+m)})-F_n(X_{(i-m)})}&\leq \frac{|F_n(X_{(i+m)})-F(X_{(i+m)})|}{F_n(X_{(i+m)})-F_n(X_{(i-m)})}
+\frac{|F_n(X_{(i-m)})-F(X_{(i-m)})|}{F_n(X_{(i+m)})-F_n(X_{(i-m)})}\nonumber\\&
+\frac{F_n(X_{(i+m)})-F_n(X_{(i-m)})}{F_n(X_{(i+m)})-F_n(X_{(i-m)})}\leq \frac{n}{2m}\sup_x|F_n(x)-F(x)|\nonumber\\&+\frac{n}{2m}\sup_x|F_n(x)-F(x)|+1\nonumber
\\&\leq\frac{n}{m}O\left(\sqrt{\frac{\log\log n}{n}}\right)+1,
\end{align} 
Now, since $\lim_{n\to\infty}\frac{\sqrt{n\log\log n}}{m}=0$, we can see that 
\begin{align}\label{5t}
\lim_{n\to\infty}\frac{\frac{2m}{n}}{F(X_{(i+m)})-F(X_{(i-m)})}=\lim_{n\to\infty}\frac{F_n(X_{(i+m)})-F_n(X_{(i-m)})}{F(X_{(i+m)})-F(X_{(i-m)})}\geq 1.
\end{align}
\eqref{4t}, \eqref{5t} and condition of $\sup_xf(x)<M<\infty$ result that for enough large $n$ we can write
\begin{align}\label{}
J_{11}\leq\frac{M^2}{n}\sum_{i=1}^{n}\left[\frac{(\frac{2m}{n})^2}{(F(X_{(i+m)})-F(X_{(i-m)}))^2}-1\right].
\end{align}
Now, we prove that when $m\to \infty$,  
\begin{align}
\lim_{n\to\infty}E\left(\frac{M^2}{n}\sum_{i=1}^n\left[\frac{(\frac{2m}{n})^2}{(F(X_{(i+m)})-F(X_{(i-m)}))^2}-1\right]\right)=0.
\end{align}
Note that
\begin{align}
\lim_{n\to\infty}&E\left(\frac{M^2}{n}\sum_{i=1}^n\left[\frac{(\frac{2m}{n})^2}{(F(X_{(i+m)})-F(X_{(i-m)}))^2}-1\right]\right)=
\lim_{n\to\infty}E\left(\frac{M^2}{n}\sum_{i=1}^m\left[\frac{(\frac{2m}{n})^2}{(F(X_{(i+m)})-F(X_{(1)}))^2}-1\right]\right)\nonumber\\&+
\lim_{n\to\infty}E\left(\frac{M^2}{n}\sum_{i=m+1}^{n-m}\left[\frac{(\frac{2m}{n})^2}{(F(X_{(i+m)})-F(X_{(i-m)}))^2}-1\right]\right)\nonumber\\&+
\lim_{n\to\infty}E\left(\frac{M^2}{n}\sum_{i=n-m+1}^n\left[\frac{(\frac{2m}{n})^2}{(F(X_{(n)})-F(X_{(i-m)}))^2}-1\right]\right).
\end{align}
Next notice that $F(X_{(1)})$, $\ldots$, $F(X_{(n)})$ is an orderd sample from uniform distribution and the random variable $F(X_{(i+j)})-F(X_{(i)})$ has the beta distribution with parameters $j$ and $n-j+1$, so we have
\begin{align}
\lim_{n\to\infty}E\left(\frac{M^2}{n}\sum_{i=1}^n\left[\frac{(\frac{2m}{n})^2}{(F(X_{(i+m)})-F(X_{(i-m)}))^2}-1\right]\right)\nonumber
\\=\lim_{n\to\infty}M^2\Big(\frac{2m^2(n-1)(2m^2+mn+2m-2n}{n^2(m-1)(m-2)(2m-1)}-1\Big)=0.
\end{align}
Since convergence in expectation to 0
implies convergence in probability to 0 for nonnegative random variables, thus
\begin{align}\label{6t}
J_{11}\stackrel{p}\to 0.
\end{align}
Note that by \eqref{3t} and the law of large numbers (LLN), we can write: 
\begin{align}\label{7t}
J_{12}=\frac{1}{n}\sum_{i=1}^nf^2(X_i)=\int f^2(x)dF_n(x)\stackrel{p}\to\int f^3(x)dx.
\end{align}
\eqref{2t}, \eqref{6t} and \eqref{7t} conclude that 
\begin{align}\label{8t}
J_1\stackrel{p}\to \int f^3(x)dx.
\end{align} 
To deal with $J_2$, we notice that since
$\xi_{im}\leq \frac{\frac{2m}{n}}{X_{(i+m)}-X_{(i-m)}}$, we can write
\begin{align*}
|J_2|&\leq\frac{2}{n}\sum_{i=1}^n\frac{(\frac{2m}{n})^2}{(X_{(i+m)}-X_{(i-m)})^2}(1-\frac{c_i}{2})
\\&=\frac{2}{n}\sum_{i=1}^n\left(\frac{2m}{n}\right)^2\left(\frac{F(X_{(i+m)})-F(X_{(i-m)}}{X_{(i+m)}-X_{(i-m)}}\right)^2
\frac{1-\frac{c_i}{2}}{(F(X_{(i+m)})-F(X_{(i-m)})^2}.
\end{align*}
From \eqref{3t} and boundedness of $f(X_i)$, we have
\begin{align}
|J_2|&\leq\frac{2M^2}{n}\sum_{i=1}^{n}\frac{4m^2}{n^2}\frac{1-\frac{c_i}{2}}{(F(X_{(i+m)})-F(X_{(i-m)})^2}\nonumber
\\&=\frac{4M^2m}{n^3}\left[\sum_{i=1}^{m}\frac{m-i+1}{(F(X_{(i+m)})-F(X_{(1)})^2}+\sum_{i=n-m+1}^{n}\frac{m-n+i}{(F(X_{(n)})-F(X_{(i-m)})^2}\right],
\end{align} 
thus 
\begin{align}
E(J_2)&\leq \frac{8M^2m}{n^2}\sum_{i=1}^m\frac{(n-1)(m-i+1)}{(i+m-2)(i+m-3)}\nonumber\\&=\frac{8M^2m(n-1)}{n^2}\left\{\frac{m(2m-1)}{2(m-1)(m-2)}+\psi(m-2)-\psi(2m-2)\right\},
\end{align}
where $\psi(x)=\frac{\Gamma^{\prime}(x)}{\Gamma(x)}$ and  converges to zero with $n\to\infty$, $m\to\infty$ and $m/n\to 0$. Thus, $J_2$ form a series of non-positive variables with expectations approaching zero, and consequently 
\begin{align}\label{11t}
J_2\stackrel{p}\to 0.
\end{align}
\eqref{1t}, \eqref{8t} and \eqref{11t} complete the proof.  
\end{proof} 
\begin{rem}
In Theorem \ref{thm2.5}, many m's can be found that satisfy the conditions of theorem, for example $m=O({(\log\log n)}^{\alpha}n^{\beta})$ for every $\alpha\in\mathbb{R}$ and $\frac{1}{2}<\beta<1$.
\end{rem}
\section{Simulation study}
In this section, the results of the Monte Carlo studies on mean squared errors (MSEs) and  biases of our introduced estimators are presented. We consider the gamma, uniform and exponential distributions, which are considered in many references. For each sample size, 10,000 samples were generated and MSEs and biases of estimators were computed. We used the $m=[\sqrt n+0.5]$  formula to estimate the varextropy suggested by \cite{GrWi}. They used this heuristic formula for entropy estimation. In addition, we choose the standard normal density as the kernel and its corresponding optimal value of $h_n=1.06sn^{-\frac{1}{5}}$, where $s$ is the sample standard deviation.
Tables \ref{Gam}-\ref{exp} give simulated biases and MSEs of $VJV_{m, n}$, $VJD$, $VJB_n$, $VJS_n$ and $VJQ_{m, n}$ for gamma, uniform and exponential distributed samples, respectively. It is apparent from the tables that  both the absolute bias and MSE  decrease as the sample size increases. 
The bold type in these tables indicates the varextropy estimator achieves the minimal MSE (absolute bias).  
\begingroup 
\setlength\tabcolsep{2pt}
\noindent
\begin{table}[h]
\caption{MSE and Bias of estimators to estimate the varextropy $VJ(X)$ of the gamma distribution (G(2, 1)).}\label{Gam}
\footnotesize
\centering
\begin{tabular}{lccccc}
\toprule
\multicolumn{6}{c}{MSE(Bias)}\\
\cmidrule(l){2-6} 
$n$ &{$VJV_{m, n}$} &{$VJD$} &{$VJB_n$}&{$VJS_n$}&{$VJQ_{m, n}$}\\
\hline
10 & 0.0120008(0.0198195) &{\bf 0.0000030(-0.0004866)}&0.0000048(-0.0011958)&0.0000072(-0.0023191)&0.0016533(0.0065913)  \\
20 & 0.0016044(0.0116757) & {\bf 0.0000013(-0.0005814)} & 0.0000026(-0.0013097)&0.0000064(-0.0024858)&0.0000979(0.0030967)\\30 & 0.0002339(0.0074458) & {\bf 0.0000010(-0.0005538)} & 0.0000022(-0.0012534) & 0.0000062(-0.0024651)& 0.0000228(0.0020731) \\
40&0.0000904(0.0052193)&{\bf 0.0000008(-0.0005351)}& 0.0000019(-0.0012115) &0.0000060(-0.0024254)&0.0000115(0.0016434) \\
50&0.0000483(0.0039833) & {\bf 0.0000007(-0.0005132)} &0.0000017(-0.0011682) & 0.00000577(-0.0023838) &0.0000072(0.0013517) \\
 100 & 0.0000125(0.0024089) & {\bf 0.0000004(-0.0004324)} & 0.0000013(-0.0010154) &0.0000048(-0.0021808)&0.0000030(0.0010021)\\
\bottomrule
\end{tabular}
\end{table}
\endgroup
\begingroup 
\setlength\tabcolsep{2pt}
\noindent
\begin{table}[h]
\caption{MSE and Bias of estimators to estimate the varextropy $VJ(X)$ of the uniform distribution(U(0, 1))}\label{Uni}
\footnotesize
\centering
\begin{tabular}{lccccc}
\toprule
\multicolumn{6}{c}{MSE(Bias)}\\
\cmidrule(l){2-6} 
$n$ &{$VJV_{m, n}$} &{$VJD$} &{$VJB_n$}&{$VJS_n$}&{$VJQ_{m, n}$}\\
\hline
10 &   5.634314(0.3022407) &0.0019133(0.0405239)&{\bf 0.0007618(0.0191408)}& 0.0007793(0.0197080)&0.5714468(0.1247980)  \\
20 &       0.2218065(0.1563245)        & 0.0013200(0.0349457) & 0.0003397(0.0141829)&{\bf 0.0001894(0.0109698)} &0.0163038(0.0553850)\\

30 & 0.0295767(0.0993167) & 0.0010754(0.0320101) & 0.0002256(0.0119855) & {\bf 0.0001103(0.0086761)} & 0.0040972(0.0372555) \\
40&0.0113412(0.0710308)&0.0009404(0.0301090)& 0.0001722(0.0107598) &{\bf 0.0000784(0.0074447)}&0.0023016(0.0282432) \\
50&0.0058115(0.0547022) & 0.0008532(0.0287799) &0.0001409(0.0099317) &{\bf 0.0000610(0.0066770)} &0.0010221(0.0222844) \\ 
100 & 0.0015858(0.0336733) & 0.0006416(0.0251320) & 0.0000811(0.0079505) &{\bf 0.0000300(0.0048685)} &0.0003130(0.0142730)\\
\bottomrule
\end{tabular}
\end{table}
\endgroup
\begingroup 
\setlength\tabcolsep{2pt}
\noindent
\begin{table}[h]
\caption{MSE and Bias of estimators to estimate the varextropy $VJ(X)$ of the exponential distribution with
mean equal to one.}\label{exp}
\footnotesize
\centering
\begin{tabular}{lccccc}
\toprule
\multicolumn{6}{c}{MSE(Bias)}\\
\cmidrule(l){2-6} 
$n$ &{$VJV_{m, n}$} &{$VJD$} &{$VJB_n$}&{$VJS_n$}&{$VJQ_{m, n}$}\\
\hline
10 &   0.4596275(0.0899478) &{\bf 0.0002119(-0.0125190)}&0.0003007(-0.0159281)&0.0003360(-0.0167951)&0.0745941(0.0470953)  \\
20 &       0.0252052(0.0485333)        & {\bf 0.0001861(-0.0129253) }& 0.0002770(-0.0163934)&0.0003233(-0.0175661)&0.0074666(0.0248586)\\

30 & 0.0089156(0.0322723) &{\bf 0.0001754(-0.0127999)}& 0.0002678(-0.0162028) & 0.0003118(-0.0173596)& 0.0023550(0.0179326) \\
40&0.0034209(0.0234392)&{\bf 0.0001673(-0.0125957)}& 0.0002598(-0.0159925) &0.0002999(-0.0170549)&0.0012288(0.0140575) \\
50&0.0015958(0.0177894) &{\bf 0.0001614}(-0.0124132) &0.0002522(-0.0157660) & 0.0002883(-0.0167369)&0.0006660{\bf (0.0110643)} \\ 
100 & 0.0004385(0.0106603) & {\bf 0.0001415}(-0.0117113) & 0.0002277(-0.0150125) &0.0002375(-0.0151743)&0.0002418{\bf (0.0074678)}\\
\bottomrule
\end{tabular}
\end{table}
\endgroup
In the case of gamma distribution, Table  \ref{Gam}, it is observed that, the estimator $VJD$ performs best in terms of bias and MSE. Under the uniform distribution, Table \ref{Uni}, it can be seen that, the MSEs and biases of $VJS_n$ are always smaller than other proposed estimators except in $n=10$. 
According to Table \ref{exp},  under exponential distribution, estimator $VJD$ performs better than other proposed estimators in terms of MSE, 
and in terms of absolute bias, estimator $VJD$ performs better in small samples and estimator $VJQ_{m, n}$ in large samples.

Finally, based on the results of our simulation study, we can conclude that the estimators $VJD$ and $VJS_n$ have better performance.
\section{Some tests of uniformity}
In this section, we introduce a few goodness-of--fit tests of uniformity. These tests are based on our proposed varextropy estimators and their percentage points and power values are obtained by Monte Carlo simulation.
\subsection{Test statistic}
Consider the class of continuous distribution functions ${F}$ with density function $f(x)$ defined on interval $[0, 1]$. Due to the fact that the variance is always nonnegative, we have $VJ(X)=Var[-\frac{1}{2}f(X)]\geq 0$. If $f$ is a standard uniform density function, it is easy to see that $VJ(X)=0$. On the other hand, if $VJ(X)=Var[-\frac{1}{2}f(X)]=0$, then we have $f(x)=c$ on $[0, 1]$, where $c$ is a constant and from the property of density function, it is obvious that $f(x)=1$. Therefore,  we can claim that $VJ(X)$ is always nonnegative and its minimum value is uniquely obtained by the standard uniform distribution.   
The above discussion gives us the idea to use the proposed varextropy estimators as the test statistics of goodness-of-fit tests of uniformity. 

Let $X_1$, $\ldots$, $X_n$ be a random sample from a continuous distribution function $F(x)$ on $[0,1]$. The null hypothesis is $H_0: F(x)=x$ and the alternative hypothesis (denoted $H_1$) is the opposite of $H_0$. Given any significance level $\alpha$, our hypothesis-testing procedure can be defined by the critical region:  
\begin{align}
G_n=\hat VJ(X)\geq C_{1-\alpha}, 
\end{align}
where $\hat VJ(X)$ is one of the proposed estimators and $C_{1-\alpha}$ is the critical value for the test with level $\alpha$. Notice when $\hat VJ(X)\stackrel{p}{\to}VJ(X)$, under $H_0$, $G_n\stackrel{p}{\to}0$ and $H_1$, $G_n$ converges  to a number larger than zero in probability.

Due to the complexity of the calculations, it is not easy to determine the distribution of the test statistics under $H_0$, so the critical values are calculated using the Monte Carlo method. By utilizing five different varextropy estimators, we  introduce the following test statistics for testing the uniformity:  
\begin{align*}
GV_{n}&=VJV_{m, n},\\
GD_n&=VJD,\\
GB_n&=VJB_n,\\
GS_{n}&=VJS_n,\\
GQ_{n}&=VJQ_{m, n}.
\end{align*}
It is necessary to mention here that, also \cite{Alizadeh} used varentropy for testing uniformity.  They showed that $VE(X)=0$ if and only if $f(x)$ is the standard uniform density function and used their proposed varentropy estimators as the test statistics for goodness-of-fit tests of uniformity.

\subsection{Power comparisons}
In this subsection, we compare the powers of our proposed test statistics with Kolmogorov-Smirnov statistic (\cite{Kol}, \cite{Smir}), and the test statistics proposed by \cite{Alizadeh}. 
The test statistic of Kolmogorov-Smirnov is given as follows:
\begin{align}
KS=max\left(max_{1\leq i\leq n}\Big\{\frac{i}{n}-X_{(i)}\Big\}, max_{1\leq i\leq n}\Big\{X_{(i)}-\frac{i-1}{n}\Big\}\right),
\end{align}
where $X_{(1)}$,  $\ldots$, $X_{(n)}$ are order statistics. Moreover, the Alizadeh Noughabi-Shafaei Noughabi Statistics are given as:
\begin{align*}
TV_n&=\frac{1}{n}\sum_{i=1}^n\log^2\left(X_{(i+m)}-X_{(i-m)}\right)-\left[\frac{1}{n}\sum_{i=1}^n\log\left(X_{(i+m)}-X_{(i-m)}\right)\right]^2,\\
TE_n&=\frac{1}{n}\sum_{i=1}^n\log^2\left(\frac{c_im/n}{X_{(i+m)}-X_{(i-m)}}\right)-{\left[\frac{1}{n}\sum_{i=1}^n\log\left(\frac{c_im/n}{X_{(i+m)}-X_{(i-m)}}\right)\right]}^2,
\end{align*}
where $X_{(i)}=X_{(1)}$ if $i<1$ and $X_{(i)}=X_{(n)}$ if $i>n$,
\begin{align*}
TD_n=\int_{-\infty}^{+\infty}\hat f(x){[\log\hat f(x)]}^2dx-{\left[\int_{-\infty}^{+\infty}\hat f(x){\log\hat f(x)}dx\right]}^2,
\end{align*}
where $\hat f$ is the kernel density function estimation of $f$ and is defined by 
\begin{align*}
\hat f(x)=\frac{1}{nh}\sum_{j=1}^nk\left(\frac{x-X_j}{h}\right),
\end{align*}
\begin{align*}
TB_n=\frac{1}{n}\sum_{i=1}^n{(\log \hat f(X_i))}^2-{\left[\frac{1}{n}\sum_{i=1}^n{\log \hat f(X_i)}\right]}^2,
\end{align*}
where $\hat f(X_i)=\frac{1}{n}\sum_{j=1}^nk(\frac{X_i-X_j}{h})$,
\begin{align*}
TC_n=\frac{1}{n}\sum_{i=1}^n\log^2\left(
\frac{\sum_{j=i-m}^{i+m}(X_{(j)}-\bar X_{(i)})(j-i)}{n\sum_{j=i-m}^{i+m}(X_{(j)}-\bar X_{(i)})^2}
\right)-\left[\frac{1}{n}\sum_{i=1}^n\log\left(
\frac{\sum_{j=i-m}^{i+m}(X_{(j)}-\bar X_{(i)})(j-i)}{n\sum_{j=i-m}^{i+m}(X_{(j)}-\bar X_{(i)})^2}
\right)\right]^2,
\end{align*}
where $\bar X_{(i)}=\frac{1}{2m+1}\sum_{j=i-m}^{i+m}X_{(j)}$,
and 
\begin{align}
TA_n=\frac{1}{n}\sum_{i=1}^n\log^2\{\hat f(X_{(i+m)})+\hat f(X_{(i-m)})\}-{\left[\frac{1}{n}\sum_{i=1}^n\log\{\hat f(X_{(i+m)})+\hat f(X_{(i-m)})\}\right]}^2.
\end{align}
In the above statistics, $m$ is the window size as defined in Subsection 2.5. 
\begin{table}[ht] 
 \small
 \caption{Percentage points of the proposed test statistics at the level $\alpha=0.05.$}\label{Tab4}
\begin{center}
 \begin{tabular}{cccccc}
 \hline
 n&$GV_{n}$&$GD_n$&$GB_n$&$GS_n$&$GQ_{n}$\\
 \hline
 10 & 0.9102 & 0.0665&0.0574& 0.0558&0.4024\\
 20 &0.4937&0.0485&0.0374&0.0277&0.1813 \\
30&0.2845&0.0413&0.0297&0.0205&0.1365\\
40&0.1933&0.0371&0.0256&0.0169&0.0767\\
50&0.1389&0.0343&0.0226&0.0148&0.0583\\
75&0.0872&0.0296&0.0183&0.0116&0.0385\\
100&0.0715&0.0270&0.0158&0.0098&0.0320\\
\hline
\end{tabular}
\end{center}
\end{table}
\begin{table}[ht] 
\small
\caption{Power comparisons of the proposed tests at a significance level of 0.05. }\label{Tab5}
\begin{center}
\begin{tabular}{|c|c|c|c|c|c|c|}
\hline
 n& $Alternative$&$GV_{n}$&$GD_n$& $GB_n$&$GS_n$&$GQ_n$\\
 \hline
&$A_{1.5}$ & 0.0739&\bf 0.1767&0.1072& 0.1663&0.0764\\
&$A_2$ & 0.1150&\bf 0.3959&0.2033& 0.3423 &0.1327\\
&$B_{1.5}$&0.0719&0.1589&{\bf 0.1629}&0.1109&0.0456\\
10&$B_2$&0.1126&{\bf 0.3494}&0.3292&0.2167&0.0624\\
&$B_3$&0.2156&{\bf 0.7326}&0.6421&0.5038&0.1294\\
&$C_{1.5}$& 0.0808&0.0415&0.0296&0.0573&\bf 0.1182\\
&$C_2$&0.1295&0.0406&0.0300&0.0660&\bf 0.2160\\
\hline
\hline
&$A_{1.5}$&0.0884&{\bf 0.2889}&0.1575&0.2497&0.1088\\
&$A_2$&0.1684&{\bf 0.6680}&0.3648&0.5377&0.2489\\
&$B1.5$&0.0879&0.2424&{\bf 0.2543}& 0.0977&0.0444\\
20&$B_2$&0.1640&{\bf 0.5707}&0.5493&0.2611&0.1061\\
&$B_3$&0.3689&{\bf 0.9459}&0.8948&0.7059&0.3549\\
&$C_{1.5}$&0.09628&0.0373&0.0251& 0.1010&\bf 0.1812\\
&$C_2$& 0.1900&0.0440&0.0281&0.1339&\bf 0.3924\\
\hline
\hline
&$A_{1.5}$&0.1114&{\bf 0.3933}&0.2092&0.3206&0.1580\\
&$A_2$&0.2425&{\bf 0.8319}&0.5078&0.6832&0.3994\\
&$B1.5$&0.1100&0.3160&{\bf 0.3416}& 0.1080&0.0690\\
30&$B_2$&0.2333&{\bf 0.7251}&0.7132&0.3674&0.2224\\
&$B_3$&0.5500&{\bf 0.9907}&0.9736&0.8706&0.6455\\
&$C_{1.5}$&0.1236&0.0385&0.0228& 0.1409&\bf 0.2450\\
&$C_2$& 0.2827&0.0498&0.0263&0.2118&\bf 0.5757\\
\hline 
\end{tabular}
\end{center}
\end{table}
For power comparisons, we compute the powers of  our proposed tests and the powers of tests mentioned above under the following alternative distributions:
\begin{align*}
A_k&: F(x)=1-{(1-x)}^k, \hspace{0.5cm}0\leq x\leq 1\mbox\ \ {({\rm for}\ k=1.5, 2)};\\
B_k&: F(x)=\left\{ 
\begin{array}{llc}
2^{k-1}x^k,&0\leq x\leq 0.5&\\
&&\mbox\ \ {({\rm for}\ k=1.5, 2, 3)};\\
1-2^{k-1}{(1-x)}^k,& 0.5\leq x\leq 1&
\end{array}\right.
\\
C_k&: F(x)=\left\{ 
\begin{array}{llc}
0.5-2^{k-1}{(0.5-x)}^k,&0\leq x\leq 0.5&\\
&&\mbox\ \ {({\rm for}\ k=1.5, 2).}\\
0.5+2^{k-1}{(x-0.5)}^k,& 0.5\leq x\leq 1&
\end{array}\right.
\end{align*}
\cite{Step} uses these alternative distributions in his study on power comparisons of some uniformity tests. 

Based on $100,000$ repetitions, the critical values are estimated as shown in Table \ref{Tab4}.
Tables \ref{Tab5} and \ref{Tab6} contain the results of simulations of our proposed test statistics and competitive test statistics with $100,000$ repetitions for $n=10, 20, 30$ and $\alpha=0.05$. The bold type in these tables indicates the statistic achieving the maximal power. 
According to Tables \ref{Tab5} and \ref{Tab6}, the performance  of tests depends on alternative distributions.
Against alternative $A$, our proposed test $GD_{n}$ and  ${\rm KS}$-based test  give competitive powers to different samples. For sample size ($n=10$), the test based on   $GD_{n}$ is the best. For moderate sample size $(n=20)$, if alternative distribution is $A_{1.5}$, $GD_{n}$ is better and otherwise, for
$n=30$, ${\rm KS}$ test statistic performs better than the others. Alizadeh Noughabi-Shafaei Noughabi test statistics have very low power in this case. For alternative distribution $B$, our proposed test $GD_{n}$ and $GB_{n}$ are the best tests and other tests have low powers in this case, especially the performance of $TV$, $TE$, $TC$, $TA$ and $KS$ is  poor. If the alternative is $C$, the test based on $TC$ is the best test, 
while the behavior of our proposed statistic $GQ_{n}$ is better than the performance of ${\rm KS}$.  
\begin{table}[ht] 
\small
\caption{Power comparisons of the Alizadeh Noughabi-Shafaei Noughabi tests and  ${\rm KS}$ test at a  significance level of 0.05. }\label{Tab6}
\begin{center}
\begin{tabular}{|c|c|c|c|c|c|c|c|c|}
\hline
 n& $Alternative$&$TV$&$TE$& $TD$&$TB$&$TC$&$TA$&$KS$\\
 \hline
&$A_{1.5}$ & 0.0526&0.0625&0.0811&0.0814& 0.0642&0.0832&\bf 0.1542\\
&$A_2$ & 0.0828&0.1014&0.1311& 0.1244&0.0883& 0.1427&\bf 0.3871 \\
&$B_{1.5}$&0.0295&0.0213&{\bf 0.0826}& 0.0810&0.0298&0.0414&0.0367\\
10&$B_2$&0.0195& 0.0213&0.1110&{\bf 0.1226}&0.0169&0.0422&0.0447\\
&$B_3$&0.0185&0.0319&{\bf 0.1873}&0.1810&0.0264&0.0609&0.0873\\
&$C_{1.5}$&0.1627&0.1409&0.0530&0.0519&{\bf 0.1779}&0.0736&0.1148\\
&$C_2$&0.3329&0.2909&0.0733&0.0620&{\bf 0.3530}&0.0924&0.2001\\
\hline\hline
&$A_{1.5}$&0.0731& 0.1185&0.1421&0.1429&0.0804& 0.1502&\bf 0.2796\\
&$A_2$&0.1523&0.2684&0.2924&0.2992&0.1772&0.3322&\bf 0.6979\\
&$B1.5$&0.0082&0.0335&0.1526& {\bf 0.1582}&0.0196&0.0585&0.0579\\
20&$B_2$&0.0053&0.0593&0.3020&{\bf 0.3067}&0.0126&0.0848&0.1246\\
&$B_3$&0.0093&0.1317&0.5122&0{\bf .5166}&0.0385&0.1746&0.4133\\
&$C_{1.5}$&0.2930&0.2316&0.0324&0.0334&{\bf 0.2933}&0.1027&0.1501\\
&$C_2$&0.6430&0.5415&0.0425&0.0384&{\bf 0.6673}&0.1326&0.3159\\
\hline\hline
&$A_{1.5}$&0.0942& 0.1835&0.2001&0.2039&0.1155& 0.2326&\bf 0.4066\\
&$A_2$&0.2342& 0.4228&0.4429&0.4532&0.2913& 0.5028&\bf 0.8751\\
&$B1.5$&0.0053&0.0582&0.2433&\bf 0.2441&0.0192&0.1215&0.0837\\
30&$B_2$&0.0083& 0.1622&0.4851&{\bf 0.4946}&0.0339&0.2716&0.2411\\
&$B_3$&0.0299&0.3916&0.7583&{\bf 0.7629}&0.1242&0.5309&0.7524\\
&$C_{1.5}$&0.4023&0.3224&0.0259&0.0232&{\bf 0.4029}&0.0832&0.1892\\
&$C_2$&0.8224&0.7425&0.0353&0.0244&{\bf 0.8450}&0.0942&0.4383\\
\hline 
\end{tabular}
\end{center}
\end{table}
\section{Real data}
\begin{eg}
\cite{Illo} used 55 data set of babies$'$ smiling time measured in seconds. The data follows a $U(0, 23)$ distribution. We transformed this data to $U(0, 1)$ and obtained test statistics as:
\[GV_{n}=0.01735047, GD_n=0.02292441, GB_n=0.01267519, GS_n=0.002416893, GQ_{n}=0.008132604.\]
The values belong to the acceptance region. Therefore, Assumption $H_0$ is accepted, that is, the data follows the standard uniform distribution. 
\end{eg}
\begin{eg}
We consider the data of 20 wild lizards collected by the zoologists in the southwestern United States. The total length (mm) of each given as follows:

$179$, $157$, $169$, $146$, $143$, $131$, $159$, $142$, $141$, $130$, $142$, $116$, $130$, $140$, $138$, $137$, $134$, $114$, $90$, $114$.
Using common tests such as Kolmogorov Smirnov, Shapiro-Wilk and Anderson-Darling, it can be seen that the data follows a normal distribution.  
Now we want to transform the data into standard uniform distribution using the integral probability transformation. Therefore if $U_i=F_0(X_i)$ for $i=1, \ldots, n$, where $F_0(X_i)$ is distribution function of the normal distribution with estaimated mean and variance of the data.  The transform data is obtained as follows:

$0.9804$, $0.8326$, $0.9408$, $0.6620$, $0.6056$, $0.3715$, $0.8562$, $0.5864$, $0.5670$, $0.3530$, $0.5864$, $0.1419$, $0.3530$, $0.5475$, $0.5081$, $0.4884$, $0.4289$, $0.1205$, 
$0.0091$, $0.1205$. 

Next, the values of the proposed test statistics based on transformed sample are computed as: 
$GV_n= 0.1453019$, $GD_n=0.03378436$, $GB_n=0.02535139$, $GS_n=0.002548617$, $GQ_n=0.03390633$.
By comparing with the percentage points of the test statistics at the 0.05 level given in Table \ref{Tab4}, which are equal to $0.4937$, $0.0485$, $0.0374$, $0.0277$ and $0.1813$, respectively, since the values of the test statistics are smaller than the corresponding critical
values  it can be conclude that the assumption of normality of the data cannot be rejected at the 0.05 significance level.
\end{eg}
\begin{eg}
In this example, the data set symbolizes the simulated strengths of glass fibers presented by \cite{MMan}

$1.014$, $1.081$, $1.082$, $1.185$, $1.223$, $1.248$, $1.267$, $1.271$, $1.272$, $1.275$, $1.276$, $1.278$, $1.286$, $1.288$, $1.292$, $1.304$, $1.306$, $1.355$, $1.361$, $1.364$, $1.379$,
$1.409$, $1.426$, $1.459$, $1.460$, $1.476$, $1.481$, $1.484$, $1.501$, $1.506$, $1.524$, $1.526$, $1.535$, $1.541$, $1.568$, $1.579$, $1.581$, $1.591$, $1.593$, $1.602$, $1.666$, $1.670$, $1.684$, 
$1.691$, $1.704$, $1.731$, $1.735$, $1.747$, $1.748$, $1.757$, $1.800$, $1.806$, $1.867$, $1.876$, $1.878$, $1.910$, $1.916$, $1.972$, $2.012$, $2.456$, $2.592$, $3.197$, $4.121$.
\end{eg}
The values of the proposed test statistics after obtaining transformed sample are $GV_n= 0.4381156$, $GD_n=0.05200372$, $GB_n=0.03980463$, $GS_n=0.03055442$, $GQ_n=0.05855942$.
By looking at the Table \ref{Tab4}, we can see that because the test statistics are larger than the percentage points obtained for n=63, therefore, the assumption of normality will be rejected at the significance level of $\alpha=0.05$.
\begin{eg}
The following data had presented in  \cite{FFul} and represents the strength of glass for aircraft window.  

$18.83, 20.8, 21.657, 23.03, 23.23, 24.05, 24.321, 25.5, 25.52, 25.8, 26.69, 26.77, 26.78, 27.05, 27.67,$\\
$29.9, 31.11, 33.2, 33.73, 33.76, 33.89, 34.76,
35.75, 35.91, 36.98, 37.08, 37.09, 39.58, 44.045, 45.29, 45.381.$ 
\end{eg}
By using Kolmogorov-Smirnov (K-S) distance test statistic, \cite{Alsh} showed that, among the distributions fitted to the data,  including the exponential and A distributions, A distribution have the best fitted for data. Now by calculating the proposed test statistics after obtaining transformed sample, we show that A distribution, unlike the exponential distribution, is a suitable distribution to fit the data. 
\begingroup 
\setlength\tabcolsep{2pt}
\noindent
\begin{table}[h]
\caption{The CDF, MLE and the proposed test statistics.}
\footnotesize
\centering
\begin{tabular}{lccccccc}
\toprule
Model &CDF&MLE&{$GV_{n}$} &{$GD_n$} &{$GB_n$}&{$GS_n$}&{$GQ_{n}$}\\
\hline
 Exponentil & $1-\exp(-\beta x)$&
0.032

&0.7052102&0.2151081&0.1541674&0.08482313&0.07347814  
\\
 A& $\exp(\frac{1}{\beta}(1-\exp(\frac{\beta}{x})))$&
125.662

&0.1404678&0.02569515&0.0144915&0.01595046&  0.07580088\\
\hline
\end{tabular}
\end{table}
\endgroup
Since, for A distribution, test statistics are lower than critical values $0.1875, 0.04088, 0.028615, 0.020166$ and $0.07852028$, thus it 
appears to fit the data very well. 
\section{Conclusion}
In this paper, we proposed estimators of varextropy. The consistency of the proposed estimators is proved. 
We employed Monte-Carlo simulation to derive the mean squared error (MSE) and bias of the estimators across various distributions. In addition, we assessed the suggested estimators based on their bias and MSE characteristics.
We introduced some goodness-of-fit tests of uniformity using  the proposed estimators. Then, we obtained the critical and power values of the proposed tests by Monte Carlo simulation and compared the powers with the competitive powers in \cite{Alizadeh} and Kolmogorov-Smirnov. Our simulation 
findings demonstrated that
the suggested tests perform well in comparison with the other tests of uniformity. 
\section{Declarations}
{\bf Conflict of interest} The authors declare no conflict of interest. 

%


\end{document}